\documentclass[11pt,a4paper]{article}
\usepackage{amssymb}
\usepackage{graphicx} 
\usepackage{amsmath}
\usepackage{latexsym}
\usepackage{theorem}
\newtheorem{teorema}{Theorem}[section]
\newtheorem{definicion}[teorema]{Definition}
\include{amslatex}
\newtheorem{proposicion}[teorema]{Proposition}

\newtheorem{corolario}[teorema]{Corollary}
{\theorembodyfont{\rmfamily}
}
{\theorembodyfont{\rmfamily}
}

\numberwithin{equation}{section}

\include{amsmath}

\begin{document}
\begin{title}
{\LARGE {\bf On a rigidity condition for Berwald Spaces}}
\end{title}
\maketitle
\author{

\begin{center}
Ricardo Gallego Torrome\footnote{r.gallego.torrome@lancaster.ac.uk\\
Partially supported during the writing of this note by EPRSC and
Cockcroft Institute}

Department of Physics, Lancaster University,
\\
Lancaster, LA1 4YB, United Kingdom \& The Cockcroft Institute, UK\\[3pt]

and

Fernando Etayo\footnote{ etayof@unican.es. F.E is partially supported by Project MTM2008-01386/MTM (Spain).}

Departamento de Matem\'{a}ticas, Estad\'{\i}stica y
Computaci\'{o}n,
\\
Facultad de Ciencias, Universidad de
Cantabria, \\

Avda. de los Castros, s/n, 39071 Santander, Spain

\end{center}}
\bigskip
\bigskip

\textbf{Abstract.} We show that which  that for a Berwald structure, any Riemannian structure that is preserved
by the Berwald connection leaves the indicatrix invariant under horizontal parallel transport. We also obtain
the converse result: if $({\bf M},F)$ is a Finsler structure such that there exists a Riemannian structure that
leaves invariant the indicatrix under parallel transport of the associated Levi-Civita connection, then the
structure $({\bf M},F)$ is Berwald.  As application, a necessary condition for pure Landsberg spaces is
formulated. Using this criterion we provide an strategy to solve the existence or not of pure Landsberg
surfaces.\footnote{
M.S.C.: 53C60, 53C05,}

\newpage

\begin{center}
\textbf{Sobre una condici\'{o}n de rigidez de los espacios de Berwald}
\end{center}
\bigskip
\textbf{Resumen.} Se muestra que la conexi\'{o}n de Levi Civita de cualquier  m\'{e}trica   Riemanniana
af\'{\i}nmente equivalente a una estructura de Berwald deja invariante por transporte paralelo la indicatriz de
dicha estructura de Berwald. También se demuestra el resultado recíproco: Si $({\bf M},F)$ es una estructura de
Finsler y existe una estructura Riemanniana cuya conexi\'{o}n de Levi Civita deja invariante por transporte
paralelo la indicatriz de la estructura de Finsler, entonces $({\bf M},F)$ es de Berwald. Como aplicaci\'{o}n se
obtiene una condici\'{o}n necesaria para que una variedad sea de Landsberg pura. Y usando este criterio se
formula una estrategia para resolver el problema de la existencia de superficies de Landsberg puras.

\bigskip

{\Large \textbf{Introduction}}

\bigskip

A Riemannian structure on a manifold is given by a Riemannian metric. As is well known, the Levi Civita
connection is an important tool associated to the structure. A more general concept is that of a Finsler
structure (see {\it Definition 1.1} below). In genreal, one cannot define a Levi Civita type connection
associated to a Finsler structure.  Given a Finsler structure, one can define several linear connections on the
pull-back bundle $\pi^* {\bf TM}\longrightarrow {\bf N}$ determined by the Finsler function $F$ and additional
conditions, usually restrictions to the ``torsion". Cartan, Chern and Berwald's linear connections are notable
examples ([1]). Although the relevant merits of these connections, compared with affine connections, are quite
complicated objects. This is one of the reasons that make Finsler geometry specially difficult to investigated,
compared with the Riemannian case.

One step on the understanding of the structure of Finsler geometry
maybe achieved by the theory presented (or better suggested) in
reference [2], where it was introduced the averaged connection.
The averaged connection is obtained from a linear connection on
$\pi ^*{\bf TM}\rightarrow {\bf N}$ by an averaged procedure on a
suitable subset ${\bf {\Sigma}}\subset {\bf N}_x\,
=\pi^{-1}(x)\subset {\bf N}$. Usually this subset is defined to be
the indicatrix ${\bf I}_x$ over $x\in{\bf M}$. The average
connection is an affine connection on the tangent bundle $\pi:{\bf
TM}\rightarrow {\bf M}$. If we perform the averaged operation on
convex combinations of connections that have the same averaged
connection, the result is the same averaged connection. We called
this property convex invariance. The connection coefficients of
the averaged connection are equal to the average of the connection
coefficients of the original connection on ${\bf \Sigma}$.

The main purpose of this note is to prove a necessary condition
for Berwald spaces. The condition is obtained using the averaged
connection and the convex invariance mentioned above. In
particular two-dimensional spaces are considered.

In {\it section} 1 we introduce the basic notions of Finsler structures that we need. We will follow the
notation from Bao-Chern-Shen [1] Usually, the concepts of Finsler Geometry are introduced by using local
coordinates, but we show some intrinsic expressions (e.g., {\it propositions 1.4} and {\it 1.5}; see also [7]).
In {\it section 2} we recall the notion of average of a linear connection in the pull back bundle $\pi ^*{\bf
TM}\rightarrow {\bf N}$ and other results from [2]. In {\it section 3} we obtain {\it proposition} 3.6 which
states that for a Berwald structure, any Riemannian structure that is preserved by the Berwald connection leaves
the indicatrix invariant under horizontal parallel transport. We also obtain the converse result, {\it
proposition} 3.7: if $({\bf M},F)$ is a Finsler structure such that there exists a Riemannian structure that
leaves invariant the indicatrix under parallel transport of the associated Levi-Civita connection, then the
structure $({\bf M},F)$ is Berwald. We finish showing that these results together with the notion of convex
invariance, can be useful in the research of pure Landsberg spaces through {\it theorem} 2.9 and a criterion for
pure two dimensional Landsberg space is given.

\section{Basics notions on Finsler Geometry}
In this section we introduce the notions of Finsler geometry as well
the notation that we will use in this work. The main reference that
we follow is [1].

Let {\bf M} be an $n$-dimensional manifold and ${\bf TM}$ its
tangent bundle. If $\{x^i\}$ is a local coordinate system on {\bf
M}, the induced local coordinate system on ${\bf TM}$ is denoted by
$\{(x^i, y^i)\}$. This type of coordinate systems on ${\bf TM}$ are
called natural coordinate systems. The slit tangent bundle is ${\bf
N}={\bf TM}\setminus \{0\}$. Then we have,
\begin{definicion}
A Finsler structure $F$ on the manifold ${\bf M}$ is a
non-negative, real function  $F:{\bf TM}\rightarrow [0,\infty [$
such that
\begin{enumerate}
\item It is smooth in the slit tangent bundle ${\bf N}$.

\item Positive homogeneity holds: $F(x,{\lambda}y)=\lambda F(x,y)$
for every $\lambda >0$.

\item Strong convexity holds: the Hessian matrix
\begin{equation}
g_{ij}(x,y): =\frac{1}{2}\frac{{\partial}^2 F^2
(x,y)}{{\partial}y^i {\partial}y^j }
\end{equation}
is positive definite in ${\bf N}$. The fundamental and the Cartan
tensors are defined by the equations:
\begin{equation}
g(x,y):=\frac{1}{2}\frac{{\partial}^2 F^2(x,y) }{{\partial}y^i
{\partial}y^j }\, dx^i \otimes dx^j .
\end{equation}
\end{enumerate}
\end{definicion}

Given a Finsler structure $({\bf M}, F)$ it is not possible to
define a Levi-Civita connection in the general case. In order
to obtain a connection related with the structure, one has to go
to higher order bundles over {\bf M}. This is done in the following
 standard way. First, we introduce a non-linear 
connection on the bundle $\pi_N:{\bf TN}\longrightarrow {\bf N}$:
\begin{enumerate}
\item there is a splitting of each tangent space ${\bf T}_u{\bf
N}$ in
complementary sub-spaces $\mathcal{V}_u $ and 
$\mathcal{H}_u$
\begin{displaymath} {\bf T}_u {\bf N}=\mathcal{V}_u \oplus
\mathcal{H}_u ,\quad \forall \,\, u\in {\bf N}
\end{displaymath}
\item $ker_u(\pi_N)=\mathcal{V}_u,\quad \forall u\in {\bf N}$.
\end{enumerate}
This decomposition is invariant by the action of ${\bf GL}(n,{\bf
R})$, which is induced by the action of the linear group ${\bf
Gl}(n,{\bf R})$ acting freely and by the right on the tangent
bundle manifold ${\bf TM}$.

A Local basis for ${\bf T}_u {\bf N}$ is given by the distributions
\begin{displaymath}
\{ \frac{{\delta}}{{\delta} x^{1}}|_u
,...,\frac{{\delta}}{{\delta} x^{n}} |_u, \,
F\frac{\partial}{\partial y^{1}} |_u,...,F\frac{\partial}{\partial
y^{n}} |_u\},\quad \frac{{\delta}}{{\delta} x^{j}}|_u
=\frac{\partial}{\partial x^{j}}|_u
-N^{i}_{j}\frac{\partial}{\partial y^{i}}|_u ;
\end{displaymath}
where the non-linear connection coefficients $N^i_j$ must be
specified. The first $n$ elements develop the horizontal subspace
$\mathcal{H}_u$ while the second half the vertical subspace
$\mathcal{V}_u$. Similarly, for the cotangent space ${\bf T^*}_u
{\bf N}$ a dual basis is defined by
\begin{displaymath}
\{ dx^{1}|_u ,...,dx^{n}|_u, \frac{{\delta}y^{1}}{F}|_u,...,
\frac{{\delta}y^{n}}{F}|_u \},\, \frac{{\delta}y^{i}}{F}|_u
=\frac{1}{F}(dy^{i}+N^{i}_{j}dx^{j})|_u .
\end{displaymath}
The manifold $\pi^* {\bf TM}$ is a subset of the cartesian product
${\bf TM}\times{\bf N} $. One has the pull-back bundle $\pi^*{\bf
TM}\rightarrow{\bf N}$ given by the square
\begin{displaymath}
\begin{array}{ccccc}
& & \pi _{2} & & \\
& \pi ^{*} TM & \longrightarrow & TM &  \\
& & & & \\
\pi _{1} & \downarrow & & \downarrow & \pi  \\
& & & & \\
& N & \longrightarrow & M & \\
& & \pi  & &
\end{array}
\end{displaymath}
The projection on the first and second factors are
\begin{displaymath}
\pi _1:\pi ^* {\bf TM}\longrightarrow {\bf N},\quad
(u,\xi)\longrightarrow u,
\end{displaymath}
\begin{displaymath}
\pi _2 :\pi ^* {\bf TM}\longrightarrow {\bf TM},\quad
(u,\xi)\longrightarrow \xi,\quad \xi \in \pi^{-1 }_1 (u).
\end{displaymath}
Every vector field $Y$ over {\bf M} can be interpreted as a section
of the tangent bundle ${\bf TM}\rightarrow {\bf M}$ and has
associated a section $\pi^* Y$ of the vector bundle $\pi^* {\bf
TM}\rightarrow {\bf N}$. In local coordinates, the associated
$\pi^*Y$ to $Y$ is given in the following way:
\begin{displaymath}
Y=Y^i(x)\frac{\partial}{\partial x^i}\longrightarrow \pi^*
Y=Y^i(x)\pi^* \frac{\partial}{\partial x^i},\quad \pi_2
(\pi^*\frac{\partial}{\partial x^i})=\frac{\partial}{\partial x^i}.
\end{displaymath}
We also use the following lifted fundamental tensor (or fiber
metric):
\begin{equation}
\pi^* g(x,y):=\frac{1}{2}\frac{{\partial}^2 F^2(x,y)
}{{\partial}y^i {\partial}y^j }\, \pi^* dx^i \otimes \pi^* dx^j .
\end{equation}
\begin{definicion}
Let $({\bf M},F)$ be a Finsler structure. The Cartan tensor is
defined by
\begin{equation}
A(x,y):=\frac{F}{2} \frac{\partial g_{ij}}{\partial y^{k}}\,
\frac{{\delta}y^i}{F} \otimes dx^j \otimes dx^k=A_{ijk}\,
\frac{{\delta}y^i}{F} \otimes dx^j \otimes dx^k .
\end{equation}
\end{definicion}
One possible non-linear 
connection is introduced by defining the non-linear connection
coefficients as
\begin{displaymath}
\frac{N^{i}_{j}}{F}={\gamma}^{i}_{jk}\frac{y^{k}}{F}-A^{i}_{jk}
{\gamma}^{k}_{rs}\frac{y^{r}}{F}\frac{y^{s}}{F},\quad i, j, k
,r,s=1,...,n.
\end{displaymath}
The coefficients ${\gamma}^{i}_{jk}$ are defined in local
coordinates by
\begin{displaymath}
 {\gamma}^{i}_{jk}=\frac{1}{2}g^{i s}(\frac{\partial g_{sj}}{\partial
x^{k}}-\frac{\partial g_{jk}}{\partial x^{s}}+\frac{\partial
g_{sk}}{\partial x^{j}}),\quad i, j, k , s=1,...,n;
\end{displaymath}
$A^{i} _{jk}=g^{i l}A_{ljk}$ and $g^{i l}g_{l j}= \delta ^{i}
_{j}.$

As we have said, there is not a Levi-Civita connection associated to
the Finsler structure. However, there are several connections that
one can define and that play a similar role to the Levi-Civita
connection in Riemannian geometry. One of these connections is
Chern's connection, which introduced through the following theorem
([1], pg 38),
\begin{teorema}
Let $({\bf M},F)$ be a Finsler structure. The pull-back vector
bundle ${\pi}^{*}{\bf
TM}\rightarrow {\bf N}$ admits a 
unique linear connection determined by the connection 1-forms $\{
{\omega}^i _j,\,\, i,j=1,...,n \} $ such that the following
structure equations hold:
\begin{enumerate}
\item Torsion free condition,
\begin{equation}
d(dx^{i})-dx^{j}\wedge w^{i}_{j}=0,\quad i,j=1,...,n.
\end{equation}
\item Almost g-compatibility condition,
\begin{equation}
dg_{ij}-g_{kj}w^{k}_{i}-g_{ik}w^{k}_{j}=2A_{ijk}\frac{{\delta}y^{k}}{F},\quad
i,j,k=1,...,n,
\end{equation}
where $A_{ijk}$ are the components of the Cartan tensor.
\end{enumerate}
\end{teorema}
A coordinate invariant characterization of Chern's connection is
given by the following two propositions,
\begin{proposicion}
Let $({\bf M}, F)$ be a Finsler structure. Then the almost
$g$-compatibility condition of the Chern's connection is 
equivalent to the conditions
\begin{equation}
{\nabla}^{ch}_{V(\tilde{X})}\pi^*g=2A(X,\cdot,\cdot),
\end{equation}
\begin{equation}
{\nabla}^{ch}_{H(\tilde{X})} \pi^*g=0 ,
\end{equation}
where $V(\tilde{X})$ is the vertical component and $H(\tilde{X})$
the horizontal component of $\tilde{X}\in {\bf T}_u {\bf N}$.
\end{proposicion}
{\bf Proof:} we write the above equations in local coordinates,
\begin{displaymath}
\nabla^{ch} \pi^* g=d(g_{ij} (x,y))dx^i \otimes dx^j +g_{ij}
\nabla^{ch} (dx^i
)\otimes  dx^j+g_{ij} dx^i \otimes  
\nabla^{ch} dx^j=
\end{displaymath}
\begin{displaymath}
=\frac{{\partial}g_{ij}}{{\partial}y^k }{\delta} y^k \otimes dx^i
\otimes
dx^j +(-g_{lj}\omega ^l 
_{ik}-g_{li}\omega ^l _{jk} + \frac{{\delta}g_{ij}}{{\delta}x^k }
)dx^k \otimes dx^i \otimes dx^j .
\end{displaymath}
\hfill$\Box$

\begin{proposicion}
Let $({\bf M}, F)$ be a Finsler structure. The torsion-free
condition of the
Chern connection is equivalent to the 
following conditions
\begin{enumerate}
\item Null vertical covariant derivative of sections of ${\pi}^*
{\bf TM}$: let $\tilde{X} \in {\bf T}_u{\bf N}$ and $Y\in {\bf
\Gamma M}$, then
\begin{equation}
{\nabla}^{ch}_{{V(\tilde{X})}} {\pi}^* Y=0.
\end{equation}
\item  Let us consider $X,Y\in {\bf TM}$ and their horizontal
lifts $\tilde{X}$ and $\tilde{Y}$. Then
\begin{equation}
\nabla^{ch}_{\tilde{X}} {\pi}^* Y-{\nabla}^{ch}_{\tilde{Y}}{\pi}^*
X-{\pi}^* ([X,Y])=0.
\end{equation}
\end{enumerate}
\end{proposicion}
{\bf Proof:} The expression (1.10) defines a section of the bundle
$\pi^*{\bf \bf T M}$ due to the commutator term, as well as (1.9).
Therefore, it is only necessary to write the above equations in
local coordinates: the commutator term is zero when the vectors
are $X=\frac{\partial}{\partial x^i},Y=\frac{\partial}{\partial
x^j}$. Then
\begin{displaymath}
{\nabla}^{ch}_{\frac{\partial }{\partial x^i}|_u} {\pi}^*
\frac{\partial }{\partial x^j}-{\nabla}^{ch}_{\frac{\partial
}{\partial x^j}|_u}{\pi}^* \frac{\partial }{\partial
x^i}=(\Gamma^l _{ij}-\Gamma ^l _{ji})\pi^*\frac{\partial
}{\partial x^l}=0,
\end{displaymath}
because due to eq. $(1.5)$, one has $\Gamma^i _{jk}=\Gamma^i
_{kj}$ ([1], pg 39). The result follows from the characterization
of the Chern connection.\hfill$\Box$
\begin{definicion}
A Berwald space is a Finsler structure such that the coefficients
of the Chern's connection live on {\bf M}.
\end{definicion}

The non-linear connection of the Cartan type is constructed in the
following way([5]). By Finsler geodesic we mean the parameterized
curves in ${\bf M}$ that are extremal of the Finsler functional
arc-length. They are solutions of the differential equations (in
the case of unit parameterized Finslerian geodesics)
\begin{equation}
\frac{d^2  x ^i}{ds^2} +\gamma^i _{jk} (x,y)\frac{dx^k}{ds}
\frac{dx^j}{ds}=0,\,\, i,j,k=1,...,n.
\end{equation}
The associated spray coefficients are
\begin{displaymath}
G^i:=(\gamma^i _{sk} (x,y)y^k y^s),\quad i,s,k=1,...,n.
\end{displaymath}
The connection coefficients of the non-linear Cartan connection
are given by the derivative of the spray coefficients,
\begin{equation}
^CN^i _{j} =\frac{1}{2}\frac{\partial^2 }{\partial y^j}(\gamma^i
_{sk} (x,y)y^k y^s).
\end{equation}
Using this non-linear connection on ${\bf TN}\longrightarrow {\bf
N}$ we can define the linear Berwald connection through the
following propositions:
\begin{proposicion}
Let $({\bf M}, F)$ be a Finsler structure. Then the
$g$-compatibility condition of the Berwald connection is:
\begin{equation}
{\nabla}^b_{V(X)}\pi^* g=2A(X,\cdot,\cdot),
\end{equation}
\begin{equation}
{\nabla}^b_{H(X)} \pi^* g= -2\nabla^b _ {l}A(\cdot,\cdot,X),\quad
l=\frac{y^i}{F}\frac{\partial }{\partial x^i} .
\end{equation}
\end{proposicion}
\begin{proposicion}
Let $({\bf M}, F)$ be a Finsler structure. Then the Berwald
connection is torsion free:
\begin{enumerate}
\item Null vertical covariant derivative of sections of ${\pi}^*
({\bf TM})$: let $\tilde {X} \in {\bf T}_u{\bf N}$ and $Y\in {\bf
\pi^*}({\bf \Gamma M})$, then
\begin{equation}
{\nabla}^b_{{V(\tilde{X})}} {\pi}^* Y=0.
\end{equation}
\item  Let us consider $X,Y\in {\bf TM}$ and the associated vector
fields with horizontal components $X^i $ and $Y^i $, $\tilde{X}$
and $\tilde{Y}$. Then
\begin{equation}
{\nabla}^b _{\tilde{X}} {\pi}^* Y-{\nabla}^b_{\tilde{Y}}{\pi}^*
X-{\pi}^* ([X,Y])=0.
\end{equation}
\end{enumerate}
\end{proposicion}
In the case of Berwald structures, the Chern's connection and the Berwald connection coincide.
\section{Averaged Connection}
We introduced in [2] a method to obtain a linear connection over {\bf M} from the Chern connection on $\pi^*{\bf
TM}$. The structure of this averaged is natural: it is defined using only canonical maps and the given Finsler
structure. It is not unique because depend on the measured used and also the sub-manifold ${\Sigma}_x$ where we
perform the integration. Also one can see that the covariant derivative associated with the averaged connection
is the limit of a convex sum of covariant derivatives in different directions of the tangent space ${\bf
T}_x{\bf M}$. However, since they are all points on the fiver $\pi^{-1}(x)\subset {\bf N}$, the convex sum of
covariant derivatives is still covariant under the structure group $GL(2,{\bf R})$. Therefore the averaged
operation can be seen as a homomorphism between $Conn(\pi^*{\bf TM})$, the space of linear connections on
$\pi^*{\bf TM}$ and  $Conn({\bf TM})$, the space of affine connections on {\bf TM},
\begin{displaymath}
<\cdot>:Conn(\pi^*{\bf TM})\longrightarrow Conn({\bf TM})
\end{displaymath}
\begin{displaymath}
\nabla\longrightarrow <\nabla>.
\end{displaymath}
The averaged connection was introduced in ref. [2]. We review
briefly this construction. The proof can be found in the original
reference [2], although for convenience of the reader we indicate
the basics steps here too.

Let $\pi^* ,\pi _1 ,\pi _2$ be the canonical projections of the
pull-back bundle $\pi^* {\bf \Gamma M}$, being ${\bf \Gamma M}$ a
tensor bundle over ${\bf M}$:
\begin{displaymath}
\begin{array}{ccccc}
& & \pi _{2} & & \\
& \pi ^{*} {\bf \Gamma M} & \longrightarrow & {\bf \Gamma M} &  \\
& & & & \\
\pi _{1} & \downarrow & & \downarrow & \pi  \\
& & & & \\
& {\bf N} & \longrightarrow & {\bf M} & \\
& & \pi  & &
\end{array}
\end{displaymath}
$\pi^* _u {\bf \Gamma M}$ denotes the fiber over $u\in {\bf N}$ of
the bundle $\pi^* {\bf  \Gamma M}$ and ${\bf \Gamma}_x {\bf M}$
are the tensors over $x\in {\bf M}$, being $S_x \in {\bf \Gamma}_x
{\bf M}$ a generic element of ${\bf \Gamma}_x {\bf M}$. $S_u$ is
the evaluation of the section $S$ of the bundle $\pi^* {\bf \Gamma
M}$ at the point $u\in{\bf N}$. The indicatrix at the point $x\in
{\bf M}$ is the compact submanifold
\begin{displaymath}
{\bf I}_x :=\{y\in {\bf T}_x {\bf M}\,\,|\,\, F(x,y)=1\}\subset
{\bf T}_x{\bf M}.
\end{displaymath}

Let us consider the element $S_u \in \pi^* _u {\bf \Gamma} {\bf
M}$ and the tangent vector field $\tilde{X}$ of the horizontal
path $\tilde{\gamma}:[0,1]\longrightarrow  {\bf N}$ 
connecting the points $u\in {\bf I}_x$ and $v\in {\bf I}_z$. The
parallel transport of the Chern connection along
${\tilde{\gamma}}$ of a section $S\in \pi^* {\bf TM}$ is denoted
by $\tau _{{\tilde{\gamma}}}S$; the parallel transport along
${\tilde{\gamma}}$ of the point $u\in {\bf I}_x$ is by definition
$\tau _{{\tilde{\gamma}}} (u)={\tilde{\gamma}} (1)\in
\pi^{-1}(z)$; the horizontal lift of a path is defined using the
non-linear connection in ${\bf N}$.

The following is a standard result, although a simpler proof can
also be found in [2],
\begin{proposicion}
Let $({\bf M}, F)$ be a Finsler structure
and ${\tilde{\gamma}}:[0,1]\longrightarrow {\bf N}$ 
the horizontal lift of a path ${\gamma}:[0,1]\longrightarrow {\bf
M}$ joining $x$ and $z$ points in {\bf M}. Then  $F(x,y)$ is
invariant by the horizontal parallel transport of the
Chern's connection. In particular, let us consider the indicatrixes over $x$ and $z$ ${\bf 
I}_x\subset {\bf T}_x{\bf M}$ and ${\bf I}_z\subset {\bf T}_z{\bf
M}$. Then ${\tau}_{\tilde{\gamma}}\pi^* ({\bf I}_x )=\pi^* {\bf
I}_z. $ Therefore the horizontal parallel transport maps ${\bf
I}_x$ into ${\bf I}_z$ as submanifolds of {\bf N}.
\end{proposicion}
{\bf Proof:} Let $\tilde{X}$ be the horizontal lift in ${\bf TN}$
of the tangent
vector field $X$ along the path ${\gamma}\subset {\bf M}$ joining $x$ and $z$, $S_1, S_2 \in 
\pi^* ({\bf T}_x {\bf M})$. Then {\it corollary} $2.5$ implies ${\nabla}_{\tilde{X}} g (S_1 
,S_2)=2 A(\tilde{X},S_1 ,S_2)=0$ because the vector field $\tilde{X}$ is horizontal and the 
Cartan tensor is evaluated in the first argument. Therefore the value of the Finslerian norm 
$F(x,y)=\sqrt{g_{ij}(x,y)\, y^i \, y^j}$, $y\in {\bf T}_x {\bf M},\, Y$ with $Y=\pi^* y$ is 
conserved by horizontal parallel transport,
\begin{displaymath}
 {\nabla}_{\tilde{X}} (F^2(x,y))={\nabla}_{\tilde{X}} (g(x,y))(Y,Y)+2 
g(x,{\nabla}_{\tilde{X}} Y)=0,
\end{displaymath}
being $\tilde{X}\in {\bf TN}$ an horizontal vector. The first term is zero because the above 
calculation.
The second term is zero because of the definition of parallel transport of sections $\nabla 
_{\tilde{X}} Y=0$.
In particular the indicatrix ${\bf I}_x$ is mapped to ${\bf I}_z$ because parallel transport 
is a diffeomorphism. \hfill$\Box$

{\bf Remark} A similar statement also holds for the linear Cartan connection $\nabla^c$ because it is a
$g$-compatible connection. For the linear Berwald connection $\nabla^b$ the result is not true for general
Finsler structure, because it is not $g$-compatible. However, in the case of Berwald structure {\it proposition
1.8} holds for $\nabla^b$ because both Cartan and Berwald connections coincide.

Let us consider $\pi^* _v {\bf \Gamma M}$ a fiber over $v\in {\bf
N}$ and the tensor space over $x$, the fiber ${\bf \Gamma}_x {\bf
M}$. For each $S_x \in{\bf \Gamma} _x {\bf 
M}$ and $\, v\in \pi^{-1}(z)$, $z\in {\bf U}$ we consider the
isomorphisms
\begin{displaymath}
\pi _2 | _v :{{\pi^*  _v {\bf \Gamma} {\bf M}}}\longrightarrow
{\bf \Gamma} _z {{\bf  M}},\quad S_v\longrightarrow S _z
\end{displaymath}
\begin{displaymath}
\pi ^*  _v :{\bf \Gamma} _z{{\bf M}}\longrightarrow {{\pi^*  _v
{\bf {\Gamma}M}}},\quad S _z\longrightarrow  \pi^* _v S_z .
\end{displaymath}
\begin{definicion}
Let $({\bf M}, F)$ be a Finsler structure, $\pi (u)=x$ and  $f\in {\bf 
\mathcal{F}M}$. Then $\pi^* f\in \mathcal{F}(\pi^*{\bf TM})$ is
defined by
\begin{equation}
\pi^* _u f=f(x),\quad \forall u\in {\bf I}_x\,\subset \pi
^{-1}(x)\,\subset {\bf N}.
\end{equation}
\end{definicion}
Let us denote the horizontal lifted operator in the following way:
\begin{displaymath}
\iota:{\bf T}{\bf M}\longrightarrow {\bf T}{\bf N}, \quad
  X=X^i\frac{{\partial}}{{\partial}x^i}|_x \longrightarrow
\iota(X)=X^i\frac{{\delta}}{{\delta}x^i}|_{u}:=X^i \delta_i ,
\end{displaymath}
\begin{equation}
u\in {\bf I}_x\,\subset \pi ^{-1}(x)\,\subset {\bf N},
\end{equation}
and the horizontal lift, defined by the non-linear connection
${N^i _j}$,
\begin{displaymath}
\iota :{\bf TM}\rightarrow {\bf TN},\,\, \iota (X)=\tilde{X},\,|\,
\tilde{X} \in \mathcal{H},
\end{displaymath}
such that if $\rho :{\bf TN}\rightarrow {\bf N}$ is the canonical
projection, $\pi\cdot\rho({\iota}(X))=X$ for $X\neq 0$.
\begin{definicion}
Consider a family of operators $A_w:=\{A_w:\pi ^*_w {\bf TM} \longrightarrow \pi^*_w {\bf 
TM}\}$ with $w\in \pi ^{-1}(x)\subset {\bf N}$. The average of this family is another 
operator $A_x : {\bf T }_x {\bf M} \longrightarrow  {\bf T}_x{\bf M}$ with $x\in {\bf 
M}$ given by the action on arbitrary sections $S\in \Gamma {\bf
TM}$ by the point-wise formula
\begin{displaymath}
<A_w\,>\,:=<\pi_2 |_u A\,\pi^* _u>_u S_x =\frac{1}{vol({\bf I}_x)}\big(\int _{{\bf I}_x} \pi _2 |_u  
A_u  \pi^* _u\, dvol \big)S_x ,
\end{displaymath}
\begin{equation}
u\in {\pi^{-1}(x)},\, S_x \in {\bf \Gamma }_x {\bf M};
\end{equation}
The volume form on $dvol$ is the standard volume form induced from
the indicatrix ${\bf I}_x$ from the Riemannian volume of the
Riemannian structure $({\bf T}_x{\bf M}-\{0\},g_x )$.
\end{definicion}
\begin{proposicion}(Averaged connection of a Finsler connection [2])
Let $({\bf M}, F)$ be a Finsler structure and let us consider a
linear connection $\nabla$ on $\pi^*{\bf TM}$. Then, there is
defined on ${\bf M}$ a linear covariant derivative along $X$,
$\tilde{{\nabla}}_X $ characterized by the conditions:
\begin{enumerate}
\item $\forall  X\in {\bf T}_x {\bf M}$ and $Y\in  \Gamma{\bf
TM}$, the covariant derivative of $Y$ in the direction $X$ is
given by the following average operation:
\begin{equation}
\tilde{{\nabla}}_X Y =<\nabla>_X  Y:=<\pi _2|_u
{\nabla}_{{\iota}_u (X)} {\pi}^* _u Y\,>_u ,\,\, u\in {\bf
I}_x\,\subset \pi ^{-1}(x)\,\subset {\bf N},
\end{equation}
\item For every smooth function $f\in {\bf \mathcal{F}}({\bf M})$
the covariant derivative is given by the following average:
\begin{equation}
\tilde{{\nabla}}_X f =<\nabla>_X f:=<\pi _2 |_u \nabla _{\iota _u
(X)} \pi^* _u f>_u ,\, u\in {\bf I}_x\,\subset \pi
^{-1}(x)\,\subset {\bf N}.
\end{equation}
\end{enumerate}
\end{proposicion}
{\bf Proof:} There is a complete proof in ref. [2, {\it section 4}] of this fact. It consists on checking that
effectively $<\nabla>$ is a covariant derivative. Here we provide a different argument. This argument also holds
for different averages, like the one used in [6] or the one used more recently in [10].

The argument follows in the following way. Consider a convex sum of linear connections $t_1\nabla_1+...
t_p\nabla_p$ such that $t_1+...+t_p=1$; the connections are linear connections on {\bf M}. It is well known that
$t_1\nabla_1+... t_p\nabla_p$ is also a linear connection. Now, consider a compact manifold $\Sigma_x\subset
\pi^{-1}(X)\subset {\bf N}$ and a set of connections on {\bf M}, all of them labelled by points on ${\Sigma}$,
so there is a map $\Theta:{\bf M}\longrightarrow Mod({\bf TM})$ such that $\int_{\Sigma} \Theta=1$ and that
$\Theta\geq 0$. Then, using a limit argument of the convex sum of linear connections on {\bf M}, we have that
the averaged of the family of connections $\{ \nabla_u\}$ defines also a linear connection on {\bf M}. To apply
to our case this argument, we only need to specify that $\Sigma_x={\bf I}_x$ and that
$\Theta(u)=\pi_2|_u\nabla_{\iota_u} \pi^*$, where the right hand side must be understood for fixed $u\in {\bf
I}_x$ and as acting on sections of $\Gamma {\bf M}$.\hfill$\Box$
\paragraph{}
Let ${\nabla}$ be a linear connection on $\pi^*{\bf TM}$. Then the
generalized torsion operator acting on the vector fields $X,Y \in
{\bf TM}$ is
\begin{displaymath}
Tor_u({{\nabla}}):{\bf T}_x{\bf M}\times {\bf \Gamma
TM}\longrightarrow \pi^*_u {\bf TM}
\end{displaymath}
\begin{displaymath}
Tor_u({{\nabla}})(X,Y)=\nabla _{\iota_u(X)} \pi^*_u Y -\nabla
_{\iota_u  (Y)} \pi^*_u X -\pi^*_u [X,Y],\quad \forall u\in {\bf
N}.
\end{displaymath}
Since this definition is point-wised, we can define globally the
$Tor(\nabla)$ as the family of maps defined as before.
\begin{proposicion}
Let ({\bf M},F) be a Finsler structure and let us define a linear
connection $\nabla$ with $Tor(\nabla)=0$. Then the torsion
$Tor({\tilde{\nabla}})$ of the average connection is zero.
\end{proposicion}
{\bf Proof}: as with the proposition before, there is a proof in [2. {\it section 4}]; it is just a calculation.
However, one can see that the proof is rather direct from the definition of torsion and from the fact that
convex sum of linear connections define a linear connection. \hfill$\Box$
\section{A rigidity property of Berwald Spaces}
We start considering a generalization of some well known
properties of linear connections over ${\bf M}$ ([3], section 5.4)
to linear connections defined on the bundle $\pi ^* {\bf
TM}\rightarrow {\bf N}$.

Given two linear connections $\nabla_1$ and $\nabla_2$ on the
bundle $\pi ^* {\bf TM}\rightarrow {\bf N}$, the difference
operator
\begin{displaymath}
B:{\bf H N}\otimes \pi^* {\bf \Gamma TM}\rightarrow \pi^*{\bf
\Gamma TM}
\end{displaymath}
\begin{displaymath}
B(\iota_u(X), \pi^*_u Y)=\,^1\nabla_{\iota_u(X)}\pi^*_u Y
-\,^2\nabla_{\iota_u(X)}\pi^*Y,
\end{displaymath}
\begin{displaymath}
\forall u\in {\bf N}, \, X,Y\in {\bf \Gamma TM}
\end{displaymath}
is an homomorphism that holds the Leibnitz rule. It is essential
in this definition that we have to our disposition a non-linear
connection to define the horizontal lift $\iota_u{X}$.

The symmetric and skew-symmetric parts $S$ and $A$ of B are
defined in the following way
\begin{displaymath}
S_u\,:{\bf T}_x{\bf M}\times {\bf \Gamma TM}\longrightarrow
\pi^*_u {\bf TM}
\end{displaymath}
\begin{displaymath}
S_u(X,Y):=\frac{1}{2}\,(B(\iota_u X,\pi^*_u Y)+B(\iota_u Y,\pi^*_u
X)).
\end{displaymath}
\begin{displaymath}
\forall u\in \pi^{-1}(x),\quad X\in {\bf T}_x{\bf M}\quad ,Y\in
{\bf \Gamma TM}.
\end{displaymath}
The antisymmetric part $A$ is defined in a similar way,
\begin{displaymath}
A_u\,:{\bf T}_x{\bf M}\times {\bf \Gamma TM}\longrightarrow
\pi^*_u {\bf TM}
\end{displaymath}
\begin{displaymath}
A_u(X,Y):=\frac{1}{2}\,(B(\iota_u X,\pi^* Y)-B(\iota_u Y,\pi^*
X)),
\end{displaymath}
\begin{displaymath}
\forall u\in \pi^{-1}(x),\quad X\in {\bf T}_x{\bf M}\quad ,Y\in
{\bf \Gamma TM}.
\end{displaymath}
As for the torsion, one can define the symmetric and
skew-symmetric parts $S$ and $A$ as a family of operators, because
the above definitions are point-wise.

Consider to vector fields $X$ and $Y$ on {\bf M} such that
$[X,Y]=$. Then, the following relation holds:
\begin{displaymath}
2A_u(X,Y)=\nabla_{1(\iota_u (X))}\pi^*_u Y -\nabla_{2(\iota_u
(X))}\pi^*_u Y -(\nabla_{1(\iota_u (Y))}\pi^*_u X
-\nabla_{2(\iota_u (Y))}\pi^*_u X)=
\end{displaymath}
\begin{displaymath}
=Tor_u(\nabla_1)(X,Y) -Tor_u(\nabla_2)(X,Y).
\end{displaymath}
Since this relation holds point-wise for all $u\in \pi^{-1}(x)\in
\subset {\bf N}$ we can write
\begin{equation}
2A(X,Y)=Tor(\nabla_1)(X, Y)-Tor(\nabla_2)(X,Y).
\end{equation}
\begin{definicion}
Let $\nabla$ be a linear connection on the vector bundle
$\pi^*{\bf TM}\longrightarrow {\bf N}$ with connection
coefficients $\Gamma^i _{jk}$. The geodesics of $\nabla$ are the
solutions of the differential equations
\begin{equation}
\frac{d^2 x^i}{ds^2}+\Gamma^i _{jk}(x,\frac{dx}{ds})
\frac{dx^j}{ds}\frac{dx^k}{ds}=0,\quad i,j,k=1,...,n,
\end{equation}
where $\Gamma^i _{jk}$ are the connection coefficients of
$\nabla$.
\end{definicion}
This differential equation can be written as
\begin{equation}
\nabla_{\iota_u (X)}\pi^*_u X=0,\quad u=\frac{dx}{dt}
\end{equation}
being $X$ the unit tangent vector to the solution in the given
point. In order to check eq. $(2.3)$ one uses local coordinates.

The following propositions are direct generalizations of the
analogous results for affine connections over ${\bf M}$ ([3]).
\begin{proposicion}
Let us consider two linear connections $\nabla_1$ and $\nabla_2$
on the vector bundle ${\pi^*}{\bf TM}\rightarrow {\bf N}$ such
that the covariant derivative along vertical directions are zero.
Then the following conditions are equivalent:
\begin{enumerate}
\item The connections $\nabla_1$ and $\nabla_2$ have the same
geodesic curves in {\bf M}.

\item B(X,X)=0, where $B=\nabla_1 -\nabla_2$.

\item S=0.

\item B=A.
\end{enumerate}
\end{proposicion}
The proof follows the lines of ref. [3, pg 64-65] and it is
omitted here. However we should mention that the equivalence of
the first statement and the other requires that the covariant
derivative of sections along vertical directions must be zero.
This condition allows to define geodesics in the way we did, being
independent of the derivative of sections of ${\pi^*{\bf TM}}$
along vertical directions in ${\bf TN}$ and in this sense being
independent of type of lift, as soon as we have a complete
horizontal lift.
\begin{proposicion}
Let $\nabla_1$ and $\nabla_2$ be linear connections on the vector
bundle $\pi ^* {\bf TM}\rightarrow {\bf N}$ such that they have
null covariant derivative in vertical directions. Then
$\nabla_1=\nabla_2$ iff they have the same parameterized geodesics
and $Tor(\nabla_1)=Tor(\nabla_2)$.
\end{proposicion}
{\bf Proof:} If $\nabla_1$ and $\nabla_2$ have the same geodesics,
they have the same symmetric part (the geodesic flow determines
the symmetric part of a connection). If they have the same
torsion, then $A=0$.\hfill$\Box$
\paragraph{}
Let us consider the bundles $\pi^* {\bf TM}\rightarrow {\bf N}$ and the tangent bundle ${\bf TM\rightarrow {\bf
M}}$ endowed with a linear connection $\nabla$. The horizontal lift of $\nabla $ (or pull-back connection, ([8,
pg 57])) is a connection on ${\pi^* {\bf TM}}\rightarrow {\bf N}$ defined by
\begin{equation}
(\pi^* \nabla)_{\iota (X)} \pi^* S =\pi^*(\nabla_X S),\quad
\tilde{X}\in {\bf TM}.
\end{equation}
One can show, writing the geodesic equation in local coordinates,
that the parameterized geodesics of both connections $\pi^*
\nabla$ and $\nabla$ are the same,
\begin{displaymath}
(\pi^* \nabla)_{\iota_u (X)} \pi^*_u X =0\quad
\Leftrightarrow\quad \nabla_X X=0,
\end{displaymath}
because the possibly non-zero connection coefficients are the
same:
\begin{displaymath}
\nabla _{\partial_j} \partial_k =\Gamma ^i _{jk} \partial_i
\,\Rightarrow \, \pi^*\nabla _{\delta_j} \pi ^* \partial_k =\pi^*
(\Gamma ^i _{jk}
\partial_i)=(\Gamma ^i _{jk}\pi^*
\partial_i).
\end{displaymath}
\begin{proposicion}
Let $\nabla^{ch}$ be the Chern connection of a Finsler structure
$({\bf M, F)}$, $\nabla^{b}$ the linear Berwald connection and
consider the average connection $<\nabla^{ch}>$. Then
\begin{enumerate}
\item The structure is Berwald iff
$\pi^*<\nabla^{ch}>=\nabla^{ch}$.

\item If $\pi^*<\nabla^{b}>=\nabla^{b}$, the structure is Berwald.
\end{enumerate}
\end{proposicion}
{\bf Proof:} If $\pi^* <\nabla^{ch}>=\nabla^{ch}$, since the
induced horizontal connection $\pi^* <\nabla^{ch}>$ has the same
coefficients that $<\nabla^{ch}>$ and they live on ${\bf M}$, the
structure $({\bf M},F)$ is Berwald.

Let us suppose that the structure is Berwald. Then
$\pi^*<\nabla^{ch}>=\pi^*\,<1>\,\nabla^{ch} = \nabla^{ch}$. This
relation is checked writing the action of the average covariant
derivative on arbitrary vector sections.

An alternative proof of is the following. We know that
$Tor(\nabla^{ch} )=Tor(<\nabla^{ch}>)=0$. On the other hand, the
parameterized geodesics of $\pi^*<\nabla^{ch}>$ are the same than
the geodesics of $<\nabla^{ch}>$. But if the space is Berwald, the
geodesic equation of $<\nabla^{ch}>$ are the same than the
geodesic equation of $\nabla^{ch}$. From this fact it follows
$\pi^*<\nabla^{ch}>=\nabla^{ch}$.

To proof the second statement we follow a similar reasoning. If
$\pi^*<\nabla^{b}>=\nabla^{b}$, the Berwald connection lives on
{\bf M} and therefore the structure is Berwald. \hfill$\Box$
\begin{proposicion}
Let $({\bf M}, F)$ be a Finsler structure. Then there is an affine
equivalent Riemannian structure $({\bf M},h)$ iff the structure is
Berwald.
\end{proposicion}
{\bf Proof:} if there is an affine equivalence Riemannian
structure $h$ such that its Levi-Civita connection $\nabla^h$ has
the same parameterized geodesics as the linear Berwald connection
$\nabla^b$ and both connection have also null torsion, then both
connections are the same ([3], section 5.4) and since the
connection coefficients $^h\Gamma^i _{ij}$ live in {\bf M}, the
structure is Berwald. Conversely, if $({\bf M},F)$ is Berwald, its
Berwald connection is metrizable ([6]).\hfill$\Box$.
\paragraph{}
Recall that for Berwald spaces $\nabla^b =\nabla^{ch}$. Then,
\begin{proposicion}
Let $({\bf M},F)$ be a Berwald structure. Then any Riemannian $h$
on {\bf M} such that $\nabla^b \pi^* h=0$ then $\nabla^h$ leaves
invariant the indicatrix under horizontal parallel transport.
\end{proposicion}
{\bf Proof}: If the Riemannian structure is conserved by the
Berwald connection, $\nabla^b \pi^* h =0$.  This implies that
$<\nabla^b > h=0$. In addition, $<\nabla^b>$ is torsion free.
Therefore, $<\nabla^b>=\nabla^h$. If $\nabla^b$ leaves invariant
the indicatrix, also $\pi^* <\nabla^b>\,=\nabla^h$ leaves
invariant the structure. \hfill$\Box$

There is a converse of this result,
\begin{proposicion}
Let $({\bf M},F)$ be a Finsler structure. Then if there is a
Riemannian metric $h$ that leaves invariant the indicatrix under
the parallel transport pull-back of its Levi-Civita connection
$\pi^*\nabla^h$, the structure is Berwald.
\end{proposicion}
{\bf Proof}: Let us consider such Riemannian metric $h$ and the
associated Levi-Civita connection $\nabla^h$. The induced
connection $ \pi^* \nabla ^h$ is torsion free and its connection
coefficients in natural coordinates live on {\bf M}. In addition,
the averaged connection $<\pi^*\nabla^h>$ coincides with
$\nabla^h$, so $ \pi^* \nabla ^h=\pi^* <\pi^*\nabla ^h>=\nabla^b
$, the last equality because $\pi^* <\pi^*\nabla ^h>$ leaves
invariant the indicatrix and it is torsion-free. Therefore the
result follows because the connection $\pi^* <\pi^*\nabla ^h>$ has
coefficients living on {\bf M}.\hfill$\Box$
\section{A corollary on pure Landsberg spaces}
Let us consider a metric $h$ such that its parallel Riemannian
transport leaves invariant the indicatrix of the Finsler metric
$F$, following {\it proposition 2.7}. Then, let us consider the
interpolating set of metrics
\begin{displaymath}
F_t (x,y)=(1-t)F(x,y)+t\sqrt{h(x)_{ij}y^i y^j},\,\,
i,j=1,...,n,\,\, t\in [0,1]
\end{displaymath}
and their indicatrix,
\begin{displaymath}
{\bf I}_x(t):=\{F_t (x,y)=1,\, y\in {\bf T}_x {\bf M}\}.
\end{displaymath}
Since the metric $F$ is Berwald, all the above interpolating
metrics define indicatrix that are invariant under the Levi-Civita
connection of $h$.

Let us consider the hypothesis that each of these indicatrix
defines a submanifold of ${\bf T}_x{\bf M}$ of co-dimension 1 and
that they are non-intersecting. Therefore the union of indicatrix
$\{{\bf I}_x (t),\,\in [0,1]\}$ defines a submanifold of ${\bf
T}_x{\bf M}$ of codimension $0$ that is invariant under the
holonomy of the metric $h$. This conditions are interesting for us
because it help to provided a necessary criteria for pure
Landsberg spaces,
\begin{definicion}
A Finsler structure $({\bf M},F)$ is a Landsberg space if the
$hv$-curvature
$P$ is such 
that $\dot{A}_{ijk}=P^n _{ijk}=0$, where the vector field is
defined as $e_n =\frac{y}{F(y)}$. A pure Landsberg space is such
that it is Landsberg and it is not Berwald or locally Minkowski.
\end{definicion}
This definition that we take of Landsberg space is a bit unusual,
although can be obtained from the standard characterizations
straightforward. In particular, Landsberg space is such that ([1],
{\it section 3.4})
\begin{displaymath}
0=\dot{A}_{ikl}=-l^j\,P_{jikl}\,=\tilde{l}_j\,P^{j}_{ikl}:= P^{n}
_{ikl}.
\end{displaymath}
\begin{teorema}
Let $({\bf M},F)$ be a Landsberg space and suppose that the
averaged connection $<\nabla^{ch}>$ does not leave invariant any
compact submanifolds ${\bf I}_x(t)\subset {\bf T}_x{\bf M}$ of
codimension zero. Then the structure $ ({\bf M},F)$ is a pure
Landsberg space.
\end{teorema}
{\bf Proof:} suppose that the Landsberg space is Berwald. Then we
know from a theorem of Szabo that this linear Berwald connection
is metrizable ([6]). Then, there is a Riemannian connection
$\nabla^h$ that is identified with the average connection
$<\nabla^{ch}>$ and this is in contradiction with the hypothesis
of the theorem because $\pi^* \nabla^h =\pi^*
<\nabla^{ch}>=\nabla^{h}$ leaves invariant the set of indicatrix
${\bf I}_x(t),\, \forall t\in [0,1]$ as we show before, the union
defining a submanifold of codimension zero of ${\bf T}_x{\bf
M}$.\hfill$\Box$

In this theorem, the hypothesis of Landsberg metric $F$ can be substituted by a general Finsler metric.
Therefore, {\it theorem 4.2} is essentially a criterion for not being Berwald.

{\bf Application of the {\it theorem 4.2} in dimension 2.} Let us consider the set of possible holonomy groups
of affine free-torsion connections ([4]). Then we look for the holonomy groups that can leave invariant a
compact, foliated manifold of dimension $2$. The possible holonomy groups for averaged connection of pure
Landsberg spaces should be excluded from this list. In particular, Riemannian holonomies are excluded. Since the
torsion of the averaged connection is zero, the only candidates for the holonomy of the averaged connection in
dimension $2$ are of the form $T_{{\bf R}}\cdot SL(2,{\bf R})$ for real representations, where $T_{{\bf R}}$
denotes any connected Lie subgroup of {\bf R}. The second possibility is the whole general group $GL(2,{\bf
R})$. From this family of groups, $SL(2,{\bf R})$ and $GL(2,{\bf R})$ are the candidate that can supply the
additional Landsberg condition,
\begin{corolario}
Let ${\bf M},F)$ be a two-dimensional Finsler structure such that
the average connection is $<\nabla^{ch}>$. Then if the space is
pure Landsberg, the holonomy group of $<\nabla^{ch}>$ is
$SL(2,{\bf R})$ or $GL(2,{\bf R})$.
\end{corolario}
This result provides a strategy to solve the problem of the existence of pure Landsberg spaces in dimension $2$.
We hope that future research could reveal the existence of pure Landsberg spaces, following the direction of
{\it Corollary 2.10} (see ref. [9] for a suggestion of realization of this strategy).

A generalization of this strategy to higher dimensions can also be
fruitful, but additional techniques are required, due to the
growth of the possible holonomies.

\bigskip

{\Large \textbf{Acknowledgments}}

The authors want to express their gratitude to  the suggestions and the interest of the referees.

\end{document}